\documentclass[12pt,reqno]{amsart}
\usepackage{amsmath,amssymb,graphicx,amscd}
\usepackage{mathrsfs}
\usepackage{enumerate}
\usepackage{dsfont}
\usepackage{amsfonts}
\usepackage{fullpage}
\usepackage{mathptmx}    
\usepackage[latin1]{inputenc}

\newcommand{\be}{\begin{eqnarray*}}
\newcommand{\ee}{\end{eqnarray*}}
\newcommand{\R}{\mathbb{R}}

\newcommand{\q}{\mathbb{Q}}
\newcommand{\p}{\mathbb{P}}

\newcommand{\smallconstant}{\varepsilon}
\newcommand{\stime}{\tau}
\newcommand{\conststime}{L}
\newcommand{\corrR}{R}
\newcommand{\chol}{L}
\newcommand{\vf}{\Theta}
\newcommand{\vc}{\theta}
\newcommand{\vcc}{\hat{\theta}}
\newcommand{\xif}{\Phi}
\newcommand{\xic}{\phi}
\newcommand{\xicc}{\hat{\phi}}
\newcommand{\volfund}{V}
\newcommand{\volfundtilde}{\tilde{V}}
\newcommand{\volac}{\sigma}
\newcommand{\vola}{\Sigma}
\newcommand{\partition}{\pi}
\newcommand{\BSDEdrift}{\Lambda}
\newcommand{\Mvol}{\zeta}
\newcommand{\Mdrift}{\mu}
\newcommand{\densi}{\rho}
\newcommand{\Vone}{U}
\newcommand{\Vonetilde}{\tilde{U}}

\newcommand{\Lone}{L^1}
\newcommand{\Ltwo}{L^2}
\newcommand{\Linfty}{L^\infty}
\newcommand{\Htwo}{\mathcal{H}^2}
\newcommand{\Mtwo}{M^2}

\newcommand{\alphaad}{a}
\newcommand{\funM}{\Gamma}

\newtheorem{theorem}{Theorem}[section]

\newtheorem{lemma}[theorem]{Lemma}
\newtheorem{proposition}[theorem]{Proposition}
\theoremstyle{definition}
\newtheorem{definition}[theorem]{Definition}
\theoremstyle{remark}
\newtheorem*{remark}{Remark}

\theoremstyle{remark}

\numberwithin{equation}{section}
\begin{document}

\title{Liquidity Risk, Price Impacts and the
Replication Problem\thanks{This
work was supported in part by the Fonds québécois de la recherche
sur la nature et les technologies and NSF Grant DMS-0306194.}
}

\author{Alexandre F. Roch
}

\address{ ETH Zürich,
Departement Mathematik \\
8092 Zürich,
Switzerland \\
alexandre.f.roch@gmail.com     
}

\maketitle

\begin{abstract}
We extend a linear version of the liquidity risk model of Çetin et al.
(2004) to allow for price impacts. We show that the impact of a market order on prices depends on the size of the transaction and the level of liquidity. We obtain a simple
characterization of self-financing trading strategies and a
sufficient condition for no arbitrage. We consider a stochastic volatility model in which the volatility is partly correlated with the liquidity process and show that, with the use
of variance swaps, contingent claims whose payoffs depend on the
value of the asset can be approximately replicated in this setting. The
replicating costs of such payoffs are obtained from the solutions
of BSDEs with quadratic growth and analytical properties of these
solutions are investigated.
\end{abstract}

\section{Introduction}\label{intro}

In financial markets,
liquidity either refers to the ease with which financial
securities can be bought and sold or to the ability to trade
without triggering important changes in asset prices.
Liquidity becomes a risk factor when the magnitude of the impact of these phenomena
changes randomly over time.
Uncertainty regarding the level of liquidity in traded assets has been for a long time a critical issue for moderate to large traders. The cost of a given trading strategy
in real world situations can be substantially high when large quantities of financial assets are traded due to the consequential impact of trading on prices, and the limited and uncertain future supply and demand. In this paper, we construct an arbitrage-free model which relates levels of liquidity to trade impacts and quantify liquidity costs of strategies used for hedging claims contingent on the value of the traded asset.

The literature on liquidity risk is
large and can be mainly divided according to these two conceptual perspectives. In the first category of models, the price of an asset depends on the size of
 the transaction and the depth of the
order book. The second category includes those commonly known
as ``large trader'' models in which a large trader buys
and sells such large quantities of assets that his trades affect
the prices in a non-negligible way.
The purpose of this paper is to combine both approaches in a
unified framework and to study the problem of contingent claim replication.

Examples of recent papers in the first category
include Çetin, Jarrow and Protter \cite{CJP2004} and Çetin and
Rogers \cite{CR2005}.
 Rogers and Singh \cite{RS2007} give  a
microeconomic argument for a price which depends on size and this
is then reflected in the dynamics of self-financing strategies. They
 solve an optimal control problem in this context.

Bank and Baum \cite{BB2004},
Frey \cite{F1998} and Jarrow \cite{J1994} are examples of papers 
in which the impact of the large trader is a function of its current holdings.
In Alfonsi et al. \cite{ASS2007}, the
authors relate the impact of trades to the shape of the order book
and consider the problem of optimal liquidation by the large trader.
On the other hand, Ly Vath et al. \cite{LMP2007} study the problem of
optimal portfolio selection for a large trader who has a price impact
function and cost function of exponential form.

Our present model was in part inspired by the liquidity risk model
of Çetin, Jarrow and Protter \cite{CJP2004} (thereafter referred
to as the CJP model). In the CJP model, liquidity is introduced by
hypothesizing the existence of a supply curve $S(t,x)$ which gives,
at a given time $t$, the price per share to pay for a stock in
terms of the size $x$ of the trade. In such a model, the trader
observes the supply curve and acts as a price taker. In this setting, liquidity costs
essentially depend on the quadratic variation of the trading
strategy. The main drawback of this model is that liquidity risk
can essentially be avoided by approximating a given self-financing
trading strategy (s.f.t.s.) by a sequence $(X^n)_{n\geq 1}$ of continuous
s.f.t.s. with finite variation (FV) which incur no liquidity costs. The prices of options are
then unaffected by liquidity risk. This issue was cleverly dealt
with in Çetin et al. \cite{CST2006} by adding constraints on the
gamma of the hedging strategies. A liquidity premium is then
reflected in option prices.

Our approach is to combine both notions of liquidity risk by
hypothesizing the existence of a random linear supply curve  and by
studying the impact of trades on prices. One of the key observations made in this paper is that the magnitude of price impacts is directly
related to the amount of liquidity of the asset. This leads to a
simple characterization of self-financing trading strategies in
which the profit is partly affected by the level of liquidity.
The main goal of this paper is to study the effect of liquidity
risk on the replicating costs of contingent claims. We consider a stochastic volatility model in which the volatility process depends in part on the level of liquidity. We will see that variance swaps are the simplest hedging tools in this setting.

 The paper is
organized as follows. In Section \ref{setupliq}, we derive the impact
of trading on prices using simple principles and show that changes
in the price of an asset is directly affected by the amount of
liquidity. We then use these observations to propose a model defined
on the Brownian filtration and show it is arbitrage-free. A simple
characterization of self-financing strategies is derived to help
set up the replication problem. Section \ref{completeness} is
devoted to the main result of this paper, the replication of contingent claims using variance
swaps and the characterization of replication costs in terms of
backward stochastic differential equations with quadratic growth.
Section \ref{PropertiesSection} presents useful analytical
properties of these solutions.

\section{The Setup}\label{setupliq}

We consider an economy consisting of a risky asset (typically a
stock) which is traded through a limit order book, its associated
contingent claims and a risk-free asset.  We take the point of
view of a hedger who observes the limit order book of the stock
and makes market orders (also known as marketable limit orders).
We start by describing the supply curve the hedger would expect
to observe if he did not trade. We call it the \emph{unaffected supply curve} and denote it by $S$. It represents the limit order book that
results from \emph{all other traders' limit and market orders}. It is a
conceptual construction which is not directly observed. We will assume that the hedger's trades have a lasting impact on prices which will be added to $S$ to obtain the \emph{actual observed
supply curve}, which we denote by $S^0$.

We are given a fixed maturity $T$ and $(\Omega, \mathcal{F},
(\mathcal{F}_{t})_{0\leq t\leq T}, \p)$ a filtered
probability space
 which satisfies the usual conditions. We assume the interest rate is constant, and for simplicity we work with the discounted price processes.
The (discounted) unaffected price process is an exogenously given adapted
continuous process $S=(S_{t}(x))_{t\geq 0, x \in \mathbb{R}}$ (or
sometimes written $S(t,x)$ for convenience). $S_t(x)$ is the
 price per share for a transaction of size $x$ at time
$t$ that would be observed if the hedger did not trade before time $t$. The actual (discounted) quoted
 price per share that all market participants obtain for a transaction of size $x$ at time
$t$ is denoted by $S^0_t(x).$ We start by assuming that the \emph{unaffected supply curve} has the following linear
structure:
\begin{eqnarray}
\label{linearstruct} S_{t}(x) & = & S_{t}+M_{t} x,\mbox{ for $x
\in \mathbb{R}$}
 \end{eqnarray}
 where $(M_{t})_{t\geq0}$ and $(S_{t})_{t\geq0}$  are positive semimartingales. Note that the fact that this function is
 continuous at $x=0$ implies there is no bid-ask spread. While it is theoretically possible for $S_t(x)$ to be
 negative for some values of $x$, it is unlikely to happen in practice since the value of $M_t$ is small.
 We assume there is a measure $\q$, equivalent to
 $\p$, under which the unaffected price process $S$ is local
 martingale. As in the classical theory, this assumption will be
 sufficient to rule out arbitrage opportunities. See Theorem
 \ref{THnoarb} below in this regard.

 The assumption that the supply curve is linear is supported by the empirical study of Blais \cite{B2006} for stocks that are frequently traded in large volumes.  The study was based on a large data set of stocks traded on several different stock exchanges in the year 2003. See also Blais and Protter \cite{BP2009}.

Before we specify the precise model for $S$ and $M$ on which we
will focus, we start by detailing general characteristics that a
liquidity risk model which include price impacts should reflect.

  Equation (\ref{linearstruct}) gives us
a way to describe the limit order book. We represent it by a
density function $\densi_t(z)$ which denotes the density of the
number of shares being offered at price $z$ at time $t$, i.e.
$\int_{z_1}^{z_2} \densi_t(z) d z$ is the number of shares
offered between prices $z_1$ and $z_2$. If a trader wants to buy
$x$ shares at time $t$ through a market order then the price to pay is
$\int_{S_t}^{z_x} z \densi_t(z) d z$ in which
$z_x$ solves $\int_{S_t}^{z_x} \densi_t(z) d z=x$.
It is clear from the linear structure of the supply curve that for
any $t \leq T$ the density equals $\densi_t(z) =
    \frac{1}{2 M_t}.$ In that case, $z_x= S_t + 2 M_t x$ and the
dollar outlay for $x$ shares is $$\frac{1}{2 M_t}
\int_{S_t}^{S_t+2 M_t x} z d z = S_t x + M_t x^2 = x
S_t(x).$$ Since $\densi$ is a measure of liquidity, we can think of
$M$ as a measure of illiquidity. Indeed, the larger is $M_t$, the
higher is the liquidity cost.

We let $X_{t}$ denote the number of shares owned by the hedger at
time $t$ and $S_{t}^{0}(x)$ denote the actual asset price per share observed
in the market, which includes the impact of the hedger's
trading strategy, i.e. $S_{t}^{0}(x)$ implicitly depends on $X$. We define $S^0_t = S^0_t(0)$ as the observed quoted price.

We now describe the impact that an arbitrary market order has on the limit order book. We will then use these observations to justify our specification of $S$ and $S^0$.
First, one should
observe that if $\Delta X_t$ shares are bought at time $t$ by a
trader through a market order, then the corresponding part of the order book is used up.
This would mean that immediately after the trade the limit order
book would have a density of $0$ for prices between $S^0_t$ and
$S^0_t+2 M_t \Delta X_t$ and $\densi_t$ elsewhere since the lowest
ask price would then be $S^0_t+2 M_t \Delta X_t$ whereas the
highest bid would remain the same. In this perspective, one can see that an implicit assumption made in
the liquidity model of Çetin et al. \cite{CJP2004} is that new
limit orders to sell are placed immediately after a trade, thereby
filling up the limit order book to its previous levels since it
is assumed that trades have no impact on the supply curve. The new
observed quoted price is the same as before and the impact on
prices is non-existent in \cite{CJP2004}. Although it is reasonable to assume that
the limit order book fills up to its previous level after a trade,
it is not clear whether the gap should be filled by bid or ask
orders. For example, if the gap is filled entirely by bid orders
after the purchase of $\Delta X_t$, then the new quoted price is
shifted upwards to $S^0_t+2 M_t \Delta X_t$. In this case, the
outcome is a full impact on prices.

The empirical findings of Weber and Rosenow \cite{WR2005} showed
that in practice the impact of trading on prices is important but
can be less than the full impact described in the previous
paragraph. In fact, they showed a negative correlation between
returns and the volume of incoming limit orders which suggests
that traders respond to buying market orders by adding new limit
orders in the opposite direction. We model this phenomenon by introducing a parameter $\lambda \in [0,1]$
measuring the proportion of new bid orders (resp. ask orders)
filling up the limit order book when a trade to buy (resp. sell)
is made at time $t$. In effect, the effective impact on prices of a
trade of size $\Delta X_t$ is to shift the quoted price to
$S^0_t + 2 \lambda M_t \Delta X_t$, whereas the density level of
the order book is unaffected.

 We have
to be careful how we define the observed price process in this
setting. Indeed, when the hedger makes a trade at time $t$
the price he pays is unaffected by the impact of this current
trade whereas prices right after $t$ will be. In this sense,
$S_{t}^{0}$ will not be càdlàg in general, although
$S_{t+}^{0}$ is and includes the impact of a trade at time
$t$.

Suppose $X$ is a simple trading strategy of the form
$X=\sum_{k=0}^{k_n} \Delta_k^n X
\mathbf{1}_{[\stime_k^n,\infty)}$ in which $\Delta_{k}^{n}X=X_{\stime_{k}^{n}}-X_{\stime_{k-1}^{n}}$ for $k=1,
\dots, n$ and $\Delta_{0}^{n}X = X_0$. Then, the
observed quoted price should satisfy
$$S^0_{t}= S_t + 2 \sum_{i=0}^{k-1} \lambda  M_{\stime_i^n} \Delta_i^n X =  S_t +  2 \sum_{i=0}^{k-1} \lambda  M_{\stime_{i-1}^n} \Delta_i^n X +  2 \sum_{i=0}^{k-1} \lambda  (\Delta_{i}^n M) (\Delta_i^n X )$$
 for any $t \in (\stime_{k-1}^n,\stime_{k}^n]$. Note that the sum in the previous equation only goes up to $k-1$ since $S^0_{\stime_{k}^n}$, which represents the price per share for a trade of size 0, is not yet impacted by the trade at time $\stime_{k}^n$. The right-limit version of this process is then given by
 \begin{eqnarray}\label{rightlimit}
 S_{t+}^0  =  S_t +  2 \sum_{i=0}^{k-1} \lambda  M_{\stime_{i-1}^n} \Delta_i^n X +  2 \sum_{i=0}^{k-1} \lambda ( \Delta_{i}^n M) (\Delta_i^n X)
 \end{eqnarray} for any $t \in [\stime_{k-1}^n,\stime_{k}^n)$ when $S$ is right-continuous.
Following these observations, we define
 \begin{eqnarray}\label{tradeimpacteq}
 S^0_{t+} &=& S_t + 2 \lambda \int_0^t M_{u-}  d X_u + 2 \lambda \int_0^t d [M,X]_u
 \end{eqnarray} for all $t\leq T$, for a general semimartingale $X$. Furthermore, we  define the observed quoted price by $S^0_t =
 \lim_{s \uparrow t} S^0_{s+}.$ By assuming that the level of liquidity $\densi_t$ is unaffected by trades, we readily obtain that the supply curve is given by
 \begin{eqnarray}\label{linearstructwithimpact}
S^0_t(x) = S^0_t + M_t x
 \end{eqnarray}   for all $0<t\leq T$ and $x \in \mathbb{R}.$
 We think of $1-\lambda$ as the fraction
 of the order book which is renewed after a market order so that in
 practice the actual impact on prices is $\lambda$ times the
 full impact.

Equation \ref{tradeimpacteq} gives us a new understanding of the
causes of volatility and its relation to illiquidity. As mentioned
earlier, $S$  is
 the price process which results from limit and market orders
 of all the other market participants. The equation suggests that the impact of the market orders of each market participant is proportional to the value of $M$. The volatility of
$S$ can then be expected to be correlated in part to $M$. (Another component of the volatility of $S$ would be related to the volatility of limit orders.) In fact, many empirical works have shown that the level of liquidity is an important determinant of  the variance of log-returns. The reader is referred to the works of Farmer et al. \cite{FGLMS2004} and Weber and Rosenow \cite{WR2006} for a more detailed discussion. The observation that these authors make is that volatility is high when liquidity is low, and low when liquidity is high. Since $M$ is a measure of illiquidity, we can expect the instantaneous variance of the log-returns of the stock price to be in part correlated with $M$.  This is a key observation which will enable us to hedge
derivatives. Indeed, in the next section, we will introduce variance swaps which will be used to hedge against the liquidity risk. Since volatility is one of the most correlated quantities to liquidity risk, this is a very natural choice. See Remark \ref{remark} in this regard.

Following these observations, we consider a stochastic volatility model for $S$:
\begin{eqnarray}\label{dynMS}
 d S_t &=&  \vola_t S_t d W_{1,t},
\end{eqnarray} in which $W_1$ is a Brownian motion defined on the filtered probability
space, and $\vola_t$ is the stochastic volatility. Recall that we are working
directly under a risk neutral measure $\q$ for
unaffected prices, hence $S$ has no drift term.
We model $M$ and $\vola$ as follows.
Define $\volfund$ and $\Vone$ as the solutions of
 \begin{eqnarray}
d \Vone_t &=&  \gamma (\Vone_t + \eta) d t +  \xif(\Vone_t) d
W_{2,t},\nonumber\\
d \volfund_t &=&  \alpha (\volfund_t+a) d t +  \vf(\volfund_t) d W_{3,t}\nonumber
\end{eqnarray} in which $W = (W_{j,t})_{j\leq 3,t\leq T}$ is a
three dimensional Brownian motion defined on the filtered probability
space, and $\alpha,\gamma,\eta,a \in \mathbb{R}$.  We define $\vola^2_t =  \Vone_t + \volfund_t$ and let $M=\smallconstant \funM(\Vone)$, in which $\funM$ is strictly increasing and twice continuously differentiable.  In practice, the process $M$ takes small values compared to $\vola$, but is also an important component of the volatility process $\vola$. As a result, the constant $\smallconstant$ is typically small.

We are using a three dimensional Brownian motion since there are three different sources of risk in this model, namely the stock price, the liquidity level and the volatility, which is, in practice, only partially dependent on the level of liquidity.
The components of $W$ are typically correlated and we denote by $\corrR = \frac1t COV(W_t)$ the matrix of instantaneous correlation coefficients. We assume $\corrR$ is positive definite and we let $\chol$ be the upper triangular matrix in the Cholesky decomposition such that $\corrR^{-1}$ = $\chol^\top \chol$. We then define $B = \chol\, W$. Then $B$ is a three-dimensional Brownian motion with independent components. We denote the components of $\chol^{-1}$ by
\be
\chol^{-1} = \left(
               \begin{array}{ccc}
                 \volac_1 & \volac_2 & \volac_3 \\
                 0 & \xic_2 & \xic_3 \\
                 0 & 0 & \vc_3 \\
               \end{array}
             \right).
\ee

 We assume the functions $\vf$ and $\xif$ are chosen so that the solutions of the above stochastic differential equations are well defined. For example, one can take $\vf(v) = v^{\vcc}$ with  $\vcc = 0, \frac12$ or $1.$
Examples of stochastic
volatility models of this form are Heston \cite{H1993} ($\vcc=\frac12$), Hull and
White \cite{HW1987} ($\vcc=1$), and Detemple and Osakwe \cite{DO2000}. Other expressions for $\vola^2$ could be used, however we have chosen this particular form for its mathematical tractability and its widespread use in theory and practice.

\subsection{Self-Financing Strategies and No Arbitrage}\label{sftsNA}

In order to properly address the problem of replicating contingent
claims, we give a characterization of self-financing
strategies and establish under which condition our model is
arbitrage-free.
In our setting, the self-financing condition is as follows.
\begin{definition} Let
$\partition_{n}:t_0=\stime_{0}^{n}\leq\stime_{1}^{n}\leq\ldots\leq\stime_{k_{n}}^{n}=T$
be a sequence of random partitions tending to the identity. A pair of processes $(X_{t},Y_{t})_{ t_0 \leq t \leq T}$ is a self-financing
trading strategy (s.f.t.s.) on $[t_0,T]$  if $X$ is a semimartingale and $Y$ is an optional
process  satisfying
\begin{eqnarray} Y_{t} & = &
Y_{t_0-}-\Delta X_{t_0} S^0(t_0,\Delta X_{t_0}) -\lim_{n\rightarrow\infty}\sum_{k=1}^{k_{n}}\Delta_{k}^{n}X
S^{0}(\stime_{k}^{n},\Delta_{k}^{n}X) 1_{\{\stime_k^n \leq
t\}}\label{sfts}\end{eqnarray} where convergence is in ucp. (See
Protter \cite{P2004} for undefined terms.) Here,
$\Delta_{k}^{n}X=X_{\stime_{k}^{n}}-X_{\stime_{k-1}^{n}}$ for $k=1,
\dots, n$.

$X_t$ represents the number of shares of the asset owned by the hedger and $Y_t$ is the position in the risk-free asset  at time $t$. The interpretation is that the position in the risk-free asset at time $t$ should be equal to the position at time $t_0$ minus the cost of all the trades between $t_0$ and $t$. Here, $Y_{t_0-}$ is the value of the position in the risk-free asset before the trade at time $t_0.$
 \end{definition}
 \begin{remark}  In the classical theory, the process $X$ is
 predictable. We take $X$ in the above definition
 to be a semimartingale
 for the stochastic integral in Equation \ref{tradeimpacteq} to be well defined.
 A consequence of Proposition \ref{equivalentdef} below is that the limit in
 Equation \ref{sfts} is well-defined, and the definition of self-financing trading strategies is independent of the sequence of random partitions.
\end{remark}

Even though s.f.t.s. are defined in terms of $S^0$, they can be
characterized in terms of the exogenously given processes $M$ and
$S$ as follows:
\begin{proposition}\label{equivalentdef} Let $t_0>0$. If $(X_{s},Y_{s})_{t_0 \leq s \leq T}$ is a
self-financing trading strategy then
\begin{eqnarray} \nonumber \label{sfts2}Y_{t}+X_t (S^0_{t+} -\lambda  M_t X_t) &=&
Y_{t_0-}+X_{t_0-}( S^0_{t_0} - \lambda M_{t_0} X_{t_0-})  +\int_{t_0}^{t} X_{u-}d S_{u} \\ & & \qquad \quad -\lambda
\int_{t_0}^t X_{u-}^2 d M_u - \int_{t_0}^t
(1-\lambda) M_{u} d[X,X]_u
\end{eqnarray} for all $t_0 \leq t \leq T.$
 \end{proposition}
 \begin{proof}
 Let
$\partition_{n}:{t_0}=\stime_{0}^{n}\leq\stime_{1}^{n}\leq\ldots\leq\stime_{k_{n}}^{n}=t$
be a sequence of random partitions tending to the identity.
 The self-financing condition is
 \begin{eqnarray*}
Y_{t}
 & = & Y_{t_0-}-\Delta X_{t_0} S^0(t_0,\Delta X_{t_0}) - \lim_{n\rightarrow\infty}\sum_{i=1}^{k_{n}}\Delta_{i}^{n}X \Big(S_{\stime_{i}^{n}}^{0}
 +M_{\stime_{i}^{n}}\Delta_{i}^{n}X\Big)\end{eqnarray*} where the
 convergence is in ucp.
 We can expand the sum in the last equation to find
 \\
 $-\lim_{n\rightarrow\infty}\sum_{i=1}^{k_{n}}\Delta_{i}^{n}X\Big(S_{\stime_{i}^{n}}^{0}
 +M_{\stime_{i}^{n}}\Delta_{i}^{n}X\Big)$
 \begin{eqnarray*} & = & -\lim_{n\rightarrow\infty}\sum_{i=1}^{k_{n}}\Big(X_{\stime_{i}^{n}}S_{\stime_{i}^{n}}^{0}
 -X_{\stime_{i-1}^{n}}S_{\stime_{i-1}^{n}}^{0}\Big)
+\lim_{n\rightarrow\infty}\sum_{i=1}^{k_{n}}X_{\stime_{i-1}^{n}}\Delta_{i}^{n}S^{0}
-\lim_{n\rightarrow\infty}\sum_{i=1}^{k_{n}}M_{\stime_{i}^{n}}(\Delta_{i}^{n}X)^{2}\\
 & = & -X_t S^0_t + X_{t_0} S^0_{t_0} + \lim_{n\rightarrow\infty}\sum_{i=1}^{k_{n}}\Big(X_{\stime_{i-1}^{n}}\Delta_{i}^{n}
 S
 +2 \lambda M_{\stime_{i-1}^{n}} X_{\stime_{i-1}^{n}}\Delta_{i-1}^{n}X\Big) 
 -\lim_{n\rightarrow\infty}\sum_{i=1}^{k_{n}} M_{\stime_{i}^{n}}(\Delta_{i}^{n}X)^{2}\\
 & = & -X_t S^0_t + X_{t_0} S^0_{t_0}+2 \lambda  M_{t_0} X_{t_0} \Delta X_{t_0}- 2 \lambda M_t X_t \Delta X_t + \lim_{n\rightarrow\infty}\sum_{i=1}^{k_{n}}X_{\stime_{i-1}^{n}}\Delta_{i}^{n}S\\
 && \quad +2
 \lim_{n\rightarrow\infty}\sum_{i=1}^{k_{n}}\lambda M_{\stime_{i}^{n}}X_{\stime_{i}^{n}}\Delta_{i}^{n}X
 -  \lim_{n\rightarrow\infty}\sum_{i=1}^{k_{n}} \lambda M_{\stime_{i}^{n}}(\Delta_{i}^{n}X)^{2}
 -
\lim_{n\rightarrow\infty}\sum_{i=1}^{k_{n}}(1-\lambda)M_{\stime_{i}^{n}}(\Delta_{i}^{n}X)^{2} \\
 & = & -X_t S^0_{t+} + X_{t_0} S^0_{t_0+} +
\lim_{n\rightarrow\infty}\sum_{i=1}^{k_{n}}X_{\stime_{i-1}^{n}}\Delta_{i}^{n}S
 + \lim_{n\rightarrow\infty}\sum_{i=1}^{k_{n}}\lambda M_{\stime_{i}^{n}}\Delta_{i}^{n}X^{2}\\
&&\quad - \lim_{n\rightarrow\infty}\sum_{i=1}^{k_{n}} (1-\lambda)
M_{\stime_{i}^{n}}(\Delta_{i}^{n}X)^{2} \end{eqnarray*}
\begin{eqnarray*}
 & = & -X_t S^0_{t+} + X_{t_0} S^0_{t_0+} + \lambda M_t X_t^2 -\lambda M_{t_0} X_{t_0}^2 + \lim_{n\rightarrow\infty}\sum_{i=1}^{k_{n}}X_{\stime_{i-1}^{n}}\Delta_{i}^{n}S
 - \lim_{n\rightarrow\infty}\sum_{i=1}^{k_{n}}\lambda X_{\stime_{i-1}^{n}}^{2}\Delta_{i}^{n} M\\
 && \quad -
\lim_{n\rightarrow\infty}\sum_{i=1}^{k_{n}}(1-\lambda)
M_{\stime_{i}^{n}}(\Delta_{i}^{n}X)^{2}
\\ & = & -X_t (S^0_{t+}-\lambda M_t X_t) +X_{t_0}( S^0_{t_0+} - \lambda M_{t_0} X_{t_0}) + \int_{t_0}^{t} X_{u-}d S_{u}\\ & &\quad  - \lambda \int_{t_0}^t X_{u-}^2 d M_u- \int_{t_0}^t (1-\lambda)M_u
d[X,X]_u\end{eqnarray*} by Theorem 21 (Chapter II) of Protter
\cite{P2004} since $X$ is càdlàg.
 \end{proof}

One can think of $Y_t + x ( S^0_{t} - \lambda M_{t} x)$ as the liquidation value of a portfolio with $x$ shares at time $t$. Indeed, take $t_0=t$ and $\Delta X_t = X_{t-}$ in Equation \ref{sfts2}. Then one finds that the cash value of a position $X_{t-}$ at time $t-$ in the stock is equal to $\Delta Y_t =  X_{t-}( S^0_{t} - \lambda M_{t} X_{t-}) - (1-\lambda) M_{t} (\Delta X_t)^2$ if it is liquidated at time $t$. Furthermore, if one uses a sequence $X^n$ a continuous and FV processes converging to $X$ (this can be done by Lemma \ref{lemmacontFVapprox}), then the liquidation value converges to  $X_{t-}( S^0_{t} - \lambda M_{t} X_{t-}).$ Consequently, $\lambda M_t$ can be interpreted as the effective liquidity parameter.

Similar to the infinitely-liquid case ($M=0$), Equation
\ref{sfts2} states that the difference in the  liquidation values between time $t_0$ and $t$ 
 is equal the cumulative gains in the risky
asset $\int_{t_0}^{t}X_{u-}dS_{u},$ except that in this case there
are added costs coming from the finite liquidity of the asset.
First note that if $\lambda = 0$ we get a linear version of the
CJP model. The integral with respect to $M$ is related to the
impact of trading. If $\lambda=0$, the limit order book is automatically refilled after a market
order, as in the CJP model. At the other extreme, when $\lambda
=1$ the impact of trading is at its
fullest. It is interesting to notice that whatever the trading
strategy used an investor always has a partial benefit from the
asset becoming more liquid. Indeed, when $M_{t}$ decreases, the
associated integral is positive no matter what the sign of $X_t$
is. To understand this, it is important to remember that the hedger's trades have a permanent impact on the quoted price which is proportional to the level of liquidity. If the liquidity is low when he purchases a share and high when he sells it, the price goes up higher after his purchase then it comes down after the sale. As a result, the hedger has a partial gain from this trade. This is a typical characteristic of large trader models. Note that unless the hedger uses a trading strategy with zero quadratic variation this is only a partial benefit because there is always a liquidity cost associated to his trades.

Using Proposition \ref{equivalentdef}, for $y \in \mathbb{R}$, we define the set
$\mathcal{Z}_y$ of payoffs of maturity $T$ attainable at price $y$
by $\mathcal{F}_{T}$-measurable random variables $Y_{T}$ of the
type
\begin{eqnarray}\label{contingent} Y_{T} & = & y+
\int_{0}^{T}X_{t-}d S_{t}- \lambda \int_{0}^{T} X_{t-}^{2}d M_t -
\int_{0}^{T} (1-\lambda) M_{t} d[X,X]_{t} \nonumber\end{eqnarray}
 in which $(X_{t})_{t\geq0}$ is càdlàg with finite quadratic variation.

We will denote by
$\mathcal{Z}\stackrel{def}{=}\bigcup_{y\in\mathbb{R}}\mathcal{Z}_{y}$
 the set of all attainable payoffs. We use the following definition of
 admissibility.

\begin{definition} Let $\alphaad \geq 0$. A s.f.t.s. $(X_{t},Y_{t})_{t\geq0}$
is $\alphaad$-admissible if \begin{eqnarray*} \int_{0}^{t}X_{s-}d
S_{s}-\lambda\int_{0}^{t}X_{s-}^{2}d M_s - \int_{0}^{t}
(1-\lambda) M_{s} d[X,X]_{s} \geq-\alphaad\end{eqnarray*} for all $t \leq T$.
 The s.f.t.s. $(X_{t},Y_{t})_{\{t\geq0\}}$ is simply said to be admissible
if it is $\alphaad$-admissible for some $\alphaad\geq0$.
\end{definition}

A strategy is
 admissible if its payoff is bounded from below. In particular, this definition
 rules out doubling strategies and is well known to be a key element in the definition of arbitrage opportunities. See Delbaen and Schachermayer \cite{DS1998} in this regard.

\begin{definition} An arbitrage opportunity is an admissible s.f.t.s. whose payoff
$Y_{T}\in\mathcal{Z}_{0}$ satisfies \begin{eqnarray}
\p\{Y_{T}\geq0\}=1 & \quad\mbox{and}\quad &
\p\{Y_{T}>0\}>0.\label{noarb}\end{eqnarray}
 \end{definition}

It is already known (see \cite{CJP2004}) that the existence of a
local martingale measure for $S$ rules out arbitrage opportunities
in the CJP model. In the presence of trade impacts, the equation
for the payoff of a s.f.t.s. has an integral with respect to $M$.
Since the integrand of this integral is always negative ($-
\lambda X_{t-}^2$), then the part of the profit coming from this
integral will be negative on average if $M$ is a submartingale
under the risk neutral measure. This idea is made precise in the
following theorem which gives a sufficient condition for no
arbitrage.

\begin{theorem}\label{THnoarb} If there
exists a measure $\q \sim \p$ under
which $S$ is a $\q$-local martingale and $M$ is a
$\q$-local submartingale, then there are no arbitrage
opportunities.
\end{theorem} \begin{proof} By the Doob-Meyer
decomposition theorem there exists a $\q$-local
martingale $\widetilde{M}$ and an increasing predictable process
$A$ such that $M=\widetilde{M}+A$. Let $Z_{t}=\int_{0}^{t}X_{u-}d
S_{u}- \lambda \int_{0}^{t} X_{u-}^{2}d M_u-
 \int_{0}^{t}(1-\lambda) M_{u} d[X,X]_{u}.$ Then
$Z_{t}+\lambda \int_{0}^{t}X_{u-}^{2}dA_{u} +  \int_{0}^{t}
(1-\lambda) M_{u} d[X,X]_{u} = \int_{0}^{t}X_{u-}d S_{u}-\lambda
\int_{0}^{t}X_{u-}^{2}d\widetilde{M}_{u}\geq-\alphaad $ since $A$ is
increasing and $\int_{0}^{t} (1-\lambda) M_{u} d[X,X]_{u} \geq 0$.
Now, $S$ and $\widetilde{M}$ are $\q$-local martingales
hence $Z_{t}+\lambda \int_{0}^{t}X_{u-}^{2}dA_{u} + \int_{0}^{t}
(1-\lambda) M_{u} d[X,X]_{u}$ is also a local martingale and
because it is bounded from below it is a supermartingale.
Therefore, $Z$ is also a supermartingale and
$\mathbf{E}_{\q}Z_{T} \leq 0$. But, because
$\q\sim\p$, if $Z_{T}$ were an arbitrage
opportunity it would also satisfy Equation \ref{noarb} with
$\q$ instead of $\p$ and
$\mathbf{E}_{\q}Z_{T}>0$.
\end{proof}

In the simplest case, when $\funM(x) = x$, it suffices to take $\gamma$ and $\eta >0.$
In the case that $\funM(x) = x^2$, if $\xif(m) = \sqrt{m}$ then we need $\gamma \eta \geq \frac14;$  if $\xif(m) = m$ then we must have $\gamma \geq \frac14.$
In the remaining parts of the paper, all expectations are  with respect to
$\q.$

\section{The Replication Problem}\label{completeness}

We now turn to the problem of contingent claims replication.
 Because the presence of the processes $M$ and $\vola$ involve
 risks that cannot be hedged completely by solely trading the stock, not all payoffs are attainable when
only the underlying asset is allowed to be traded. Because these
two processes are components of the instantaneous variance of the log-returns of the stock, the
natural hedging instruments to consider are variance swaps. We
thus consider contingent claims denoted by $G_{i}$ ($i=1,2$) for
which the payoff at time $T_i>T$ ($T_1 \neq T_2$) equals the
difference between the realized variance over the time interval $[0,T_i]$ and a strike $K_i$, i.e.,
\be
G_{i,T_i}&=& \int_0^{T_i} \vola^2_s d s -K_i
= \int_0^{T_i} (\Vone_s+\volfund_s) d s -K_i.\ee To rule out arbitrage opportunities, we assume the unaffected price processes
$G^{i}$ are $\q$-martingales $(i=1,2)$, i.e.
\be
G^i_t = \mathbf{E}\left( \int_0^{T_i} \vola^2_s d s -K_i \Big| \mathcal{F}_t\right)
\ee for all $t\leq T_i.$ We further assume the $G_i$'s have a
linear supply curve, i.e. $G_{i,t}(x) = G_{i,t} + x M'_{i,t}$ for
all $x$ and $t\leq T$. Since it
is not infinitely liquid, trading $G_i$ can affect its price and
we denote by $\lambda_i M'_{i,t}$ its effective liquidity.
Typically, changes in the supply curves of the $G_i$'s will happen
less often. Hence, to keep the problem tractable, we assume that $
M'_{i,t} \equiv M'_{i,0}$ is some given positive constant for all $t
\in [0,T]$. We will see that two of these swaps are sufficient to
hedge against liquidity risks. Because we now have two more traded
assets,
$\chi_{1,t}$ denotes the number of shares of $G_1$ and $\chi_{2,t}$ the
number of shares of $G_2$ in the portfolio at time $t$. We can
easily extend the definition of s.f.t.s. to the case of three
traded securities. As shown before, s.f.t.s. $(X,\chi,Y)$ satisfy
\begin{eqnarray*}
Y_T &=& Y_t  + \int_{t}^{T}X_{u-}d S_{u} + \sum_{i=1,2} \int_t^T
\chi_{i,u-} d G_{i,u}
 - \lambda \int_{t}^{T} X_{u-}^{2}d M_{u} \\
&&  - \int_{t}^{T} (1-\lambda) M_{u} d[X,X]_u - \sum_{i=1,2}
\int_{t}^{T} (1-\lambda_i) M'_{i,u} d[\chi_{i},\chi_{i}]_u
\end{eqnarray*} for $t_0 \leq t \leq T$, when $X_{t_0} =\chi_{1,t_0}=\chi_{2,t_0}= X_T = \chi_{1,T}=\chi_{2,T}=0.$

The following proposition gives condition under which the three price processes $S, G_1, G_2$ are non-redundant. It justifies the choice of variance swaps by providing a simple explicit representation of the processes $G_i$. This result will then be used to solve the replication problem.

\begin{proposition}\label{swapsinvertible} Suppose $\alpha \neq \gamma$ and $\vf$ (resp. $\xif$) satisfies one of the following conditions:
\begin{enumerate}
  \item $\vf(v) = v^{\vcc}$ (resp. $\xif(m) = m^{\xicc}$) for $\vcc \in [0,\frac12]$ (resp. $\xicc \in [0,\frac12]$),
  \item $\vf$ (resp. $\xif$) is Lipschitz continuous.
\end{enumerate}
 Then, there exists a predictable process $\psi = (\psi_{i,j,t})_{1\leq i,j \leq 3, 0\leq t\leq T}$ in $\mathbb{R}^{2 \times 2}$ such that
\be
    G^i_t  &=& \mathbf{E} \left(
    \int_0^{T_i} \vola^2_s d s - K_i \right)  + \sum_{j=1,2,3} \int_0^t \psi_{i,j,s} d B_{j,s} \mbox{ and}\\
    S_t &=& \sum_{j=1,2,3} \int_0^t \psi_{3,j,s} d B_{j,s}
\ee   for all $t\leq T$ and $i=1,2.$ Furthermore, $(\psi_{i,j,t})_{1\leq i,j\leq 3}$ is invertible for all $t.$
\end{proposition}
\begin{proof}
Consider the process
$\volfundtilde_t:=e^{-\alpha t}(\volfund_t+a)$ for $t\leq T.$ Then,
\begin{eqnarray*}
d \volfundtilde_t & = & \sum_{i=1,2,3} e^{-\alpha t} \vc_i \vf(\volfund_t) d
B_{i,t}.
\end{eqnarray*} (We let $\vc_1=\vc_2=\xic_1=0.$) In other words, $\volfundtilde$ is a
local martingale. We first show that $\volfundtilde$ is in fact a martingale. Suppose $\vf$ is Lipschitz continuous. By the Burkholder-Davis-Gundy Inequality, there exists a positive constant $C$ such that
\be
\mathbf{E} \sup_{t\leq T} \volfundtilde_t^2 &\leq & C \mathbf{E} \int_0^T e^{-2 \alpha t} \vf(\volfund_t)^2 d t\\& \leq & C \mathbf{E} \int_0^T  \volfund_t^2 d t + C \leq  C  \int_0^T \mathbf{E} \volfund_t^2 d t + C < \infty
\ee by well known estimates of moments of solutions of stochastic differential equations with Lipschitz coefficients.  On the other hand, if $\vf(v) = v^{\vcc}$ for $\vcc \in [0,\frac12],$ then
\be
\mathbf{E} \sup_{t\leq T} \volfundtilde_t^2
  & \leq & C  \int_0^T \mathbf{E}  \volfund_t^{2 \vcc} d t
 \leq  C  \int_0^T \left(\mathbf{E} \volfund_t\right)^{2 \vcc} d t
\\& \leq &C  \int_0^T \left(e^{\alpha t} \mathbf{E} \volfundtilde_t\right)^{2 \vcc} d t
 \leq C  \int_0^T \left(e^{\alpha t} \volfundtilde_0\right)^{2 \vcc} d t < \infty
\ee since $\volfundtilde$ is a positive local martingale. Hence $\volfundtilde$ is a martingale. Similarly, we can show that the process $\Vonetilde$, defined by $\Vonetilde_t:=e^{-\gamma t}(\Vone_t+\eta)$ for $t\leq T$, is a martingale when $\xif$ satisfies Condition 1 or 2.

Suppose $\gamma, \alpha \neq 0$. Then,
\be
\mathbf{E}\left( \int_0^{T_i} \Vone_s d s | \mathcal{F}_t\right)&= & \int_0^t \Vone_s d s + \mathbf{E}\left( \int_t^{T_i} \left( e^{\gamma s} \Vonetilde_s - \eta d s \right) | \mathcal{F}_t\right)\\
&= & \int_0^t \Vone_s d s +  \int_t^{T_i}  e^{\gamma s}\left(\mathbf{E} \left( \Vonetilde_s   | \mathcal{F}_t\right)- \eta \right)d s
\\
&= & \int_0^t \Vone_s d s +  \int_t^{T_i}\left(  e^{\gamma s} \Vonetilde_t - \eta \right)d s
\\
&= & \int_0^t \Vone_s d s + \Vonetilde_t\left( \frac{e^{\gamma T_i}-e^{\gamma t}}{\gamma}\right)  - \eta(T_i - t)
\\
&= & \left( \frac{e^{\gamma T_i}-1}{\gamma}\right) \Vonetilde_0  - \eta T_i + \int_0^t\left( \frac{e^{\gamma T_i}-e^{\gamma s}}{\gamma}\right) d  \Vonetilde_s.
\ee
In particular, \be G^i_t &= &\mathbf{E}\left( \int_0^{T_i}   \Vone_s d s | \mathcal{F}_t\right) + \mathbf{E}\left( \int_0^{T_i} \volfund_s d s | \mathcal{F}_t\right) - K_i\\
&= & \left( \frac{e^{\gamma T_i}-1}{\gamma}\right)   \Vonetilde_0  - \eta   T_i + \left( \frac{e^{\alpha T_i}-1}{\alpha}\right) \volfundtilde_0  - a T_i
- K_i\\&& + \int_0^t\left( \frac{e^{\gamma T_i}-e^{\gamma s}}{\gamma}\right)   d  \Vonetilde_s + \int_0^t\left( \frac{e^{\alpha T_i}-e^{\alpha s}}{\alpha}\right) d  \volfundtilde_s \\&=& \left( \frac{e^{\gamma T_i}-1}{\gamma}\right)   \Vonetilde_0  - \eta   T_i + \left( \frac{e^{\alpha T_i}-1}{\alpha}\right) \volfundtilde_0  - a T_i
- K_i + \sum_{j=1,2,3}\int_0^t \psi_{i,j,s} d B_{j,s}\ee in which
\be\psi_{i,j,t} = \left( \frac{e^{\gamma T_i}-e^{\gamma t}}{\gamma}\right)\xic_j \xif(\Vone_t)
 +\left( \frac{e^{\alpha T_i}-e^{\alpha t}}{\alpha}\right) \vc_j \vf(\volfund_t)\ee
 for $i=1,2$ and  $j=1,2,3.$ Define
 $\psi_{3,j,t}= \volac_j \vola_t S_t$ for $j=1,2,3.$
Note that $\psi_{i,3,t} = 0$ for $i=1,2$. Since \be(\psi_{i,j,t})_{1\leq i\leq2,2\leq j\leq 3} = \left(
             \begin{array}{cc}
                \frac{e^{\gamma T_1}-e^{\gamma t}}{\gamma} &  \frac{e^{\alpha T_1}-e^{\alpha t}}{\alpha}  \\
                \frac{e^{\gamma T_2}-e^{\gamma t}}{\gamma} & \frac{e^{\alpha T_2}-e^{\alpha t}}{\alpha}  \\
             \end{array}
           \right) \left(
             \begin{array}{cc}
               \xic_2   \xif(\Vone_t) &   \xic_3   \xif(\Vone_t) \\
               0 & \vc_3 \vf(\volfund_t)  \\
             \end{array}
           \right)
\ee is invertible, so is $\psi_{t}$. In the case that $\alpha$ (resp. $\gamma$) is equal to zero, the term $\frac{e^{\alpha T_i}-e^{\alpha t}}{\alpha}$  (resp. $\frac{e^{\gamma T_i}-e^{\gamma t}}{\gamma}$) in the above matrix is replaced by $T_i-t$,  and the matrix is also invertible when $\alpha \neq \gamma.$
\end{proof}

\begin{remark}\label{remark}
The fact that the matrix $\psi$ can be explicitly obtained and shown to be invertible is the main benefit of using variance swaps to complete the market. Similar calculations for non-linear contingent claims like put and call options on the stock or the realized variance would have been much more difficult, if not impossible, to obtain. As a result, such non-linear contingent claims would make the hedging much more difficult in practice. Note that the processes $U$ and $V$ need not be martingales under the risk neutral measure, i.e. $\alpha, \gamma \neq 0$. Consequently, $\int_0^t \vola^2_s d s$ is not a martingale  and $G_{i,t} \neq \int_0^t \vola^2_s d s - K_i$ for $i=1,2$. If that were the case, one of the two variance swaps would be redundant.
\end{remark}
From now on, we assume that $\alpha\neq \gamma$ and that $\vf$ and $\xif$ satisfy the assumptions of the previous proposition.

The next lemma implies that the best way of trading is always to
use FV continuous s.f.t.s. to avoid liquidity costs coming from
the quadratic variation of $X$. In this sense, trades should
always be done at the quoted price $S(t,0)$. Note that even
though some of the liquidity costs in Equation \ref{sfts2} are
eliminated when using continuous FV strategies, liquidity risk has
not been completely eliminated from the model since the integral
$\int_t^T X_{u-}^2 d M_u$ is still present. That is the main difference between our setup and the CJP model. 

If $S$ is a special semimartingale with canonical decomposition $S=N+A$, i.e. in which $N$ is a local martingale and $A$ is a predictable and FV process, then the $\Htwo$ norm of $S$ is defined as
$$|\!|S |\!|_{\Htwo} =  |\!|\sqrt{[N]_T} |\!|_{L^2} + |\!|\int_0^T |d A_s|  |\!|_{L^2}.$$ 

\begin{lemma}\label{lemmacontFVapprox}
Let $S$ be a special semimartingale and $X$ be predictable and integrable
with respect to $S$. There exists a sequence $\{X^n\}_n$ of
bounded continuous processes with finite variation such that
$X^n_0= X^n_T=0$ and $X^n$ converges to $X$ in $\Htwo.$ In
particular, $ \int X^n d S \rightarrow \int X d S$ in
$\Htwo.$
\end{lemma}
\begin{proof}
The statement is proved in the proof of Lemma 4.1 of Çetin et al.
\cite{CJP2004}.
\end{proof}

We will see that, because of the quadratic variation term in the
equation of s.f.t.s., it is not possible to  replicate exactly in
general. Since continuous processes with finite variation have
zero quadratic variation, the previous lemma will prove to be
useful for the replication problem. Following Çetin et al. \cite{CJP2004}, we make the following
definition.
\begin{definition}\label{approxreplicationdef}
$H \in \Lone$ can be approximately replicated if there
exists a sequence $(X^{n},\chi^n,Y^{n})_{n \geq 1}$ of s.f.t.s.
such that $Y_T^{n} \rightarrow H$ in $\Lone$.
\end{definition}

In the presence of trade impacts, the process $S^0$ implicitly
depends on $X$ and its value at the maturity is \begin{eqnarray*}S_{T+}^0 = S_T+2
\lambda \int_0^T M_{u-} d X_u + 2 \lambda \int_0^T d [M,X]_u =S_T- 2
\lambda \int_0^T X_{u-} d M_u \end{eqnarray*} when $X_T=0$ and  $X_0 = 0.$ (Here
we use the time $T+$ and $X_T=0$ to make sure that the hedging strategy is
liquidated before the payoff is calculated to avoid discrepancies
between the observed asset price before and after the maturity.)
The true replication problem involves finding a s.f.t.s. $(X,Y)$
that replicates a terminal condition which itself depends on $X$.
Instead, for each $x \in \R$, we consider the replication of the terminal condition
given by $x h(\widetilde{S}_{T}^x)$ with $\widetilde{S}_{T}^x :=
S_T - 2\lambda \int_0^T x \hat{X}_{u-} d M_u$ in which $\hat{X}$ is the solution
of the replication problem in the case $\lambda=0,$ $\smallconstant = 0$ and $x=1$. Jarrow \cite{J1994} used a similar approach and
interpreted $\hat{X}_t$ as the market's perception of the option's
``delta" $X_t$. In the expression for $\widetilde{S}_{T}^x$, $x$ denotes the number of units to be replicated. Hence, the proposed approximation for the true delta for the replication of $x$ units is $x X_t.$  Proposition \ref{errorapprox} in the next section gives
an upper bound of the error introduced by this approximation. Let
us begin by giving an overview of the replication problem in this
simplified setting.

\subsection{Contingent Claims Replication Without Trade Impact and Liquidity Costs}\label{replicasubsec}

When $\lambda =
0$ and $\smallconstant=0$, the s.f.t.s. $(X_s,Y_s)_{t \leq s \leq T}$ that replicates a
payoff $H \in L^1$ satisfies
\begin{eqnarray}\label{linearBSDE}
H &=& Y_t + \int_t^T X_{u-} d S_u + \sum_{i=1,2} \int_t^T
\chi_{i,u-} d G_{i,u}.
\end{eqnarray} Also, $S^0 \equiv S$.

First, note that Equation \ref{linearBSDE} is equivalent to the
following linear backward stochastic differential equation (BSDE):
\begin{eqnarray}\label{BSDEvol}
Y_t & = & H - \sum_{j=1}^3 \int_t^T \left(\volac_j \vola_s X_s
S_s + \chi_{1,s} \psi_{1,j,s}+\chi_{2,s} \psi_{2,j,s} \right)d
B_{j,s},
\end{eqnarray}$0\leq t\leq T$. Setting
\begin{eqnarray}\label{ZXchi}Z_{j,t} =  \chi_{1,s}
\psi_{1,j,s}+\chi_{2,s} \psi_{2,j,s} + X_s \psi_{3,j,s}\end{eqnarray} for $j=1,2,3,$  the BSDE can be written as
\begin{eqnarray}\label{BSDEvol2}
Y_t & = & H - \sum_{j=1}^3 \int_t^T Z_{j,s} d B_{j,s} \quad (0\leq
t\leq T).
\end{eqnarray} When $H \in \Ltwo$, BSDE \ref{BSDEvol2} has a unique solution $(Z,Y)$
in $\Mtwo\times \Mtwo$ (see Pardoux and Peng \cite{PP1990} for
example). Since
$\psi_t$ is invertible, we can define
$X_s= \frac{Z_{1,s}}{\volac_1 \vola_s S_s}$ $(t
 \leq s \leq T)$ and $\chi_s$ by inverting Equation
\ref{ZXchi}. Then  $(X_s,\chi_s,Y_s)_{t\leq s \leq T}$ is the
solution of \ref{linearBSDE}.

\subsection{The Replication Problem With Liquidity Risk}\label{replicasubsec2}

From now on, we denote by $(\hat{X},\hat{\chi},\hat{Y})$ the solution of \ref{BSDEvol} with terminal condition $H=h(S_T).$
Recall that $\widetilde{S}_{T}^x := S_T - 2 x \lambda \int_0^T \hat{X}_u
d M_u$. The main result of this
section is the following theorem.

\begin{theorem}\label{thcomplete} Let $h : \mathbb{R}^+ \rightarrow \mathbb{R}$ be Lipschitz continuous. Then $x h(\widetilde{S}_{T}^x)$ can be
approximately replicated for all $x \in \mathbb{R}.$
\end{theorem}

\begin{proof}
Let $\conststime > 0$ and $N>0.$ Let $x \in \mathbb{R}$ and $h$ satisfy the
conditions of the theorem, and define $h^N(y) =  h(y)$ if $|y|
\leq N$ and $h^N(y)=h(N)$ otherwise. Since $h$ is continuous on
$[-N,N]$, $h^N$ is bounded. Denote this bound by $C_N$. Define
$H_T^N = x h^N(\widetilde{S}_{T}^x)$ and
\begin{eqnarray*}\stime_\conststime=\inf\{0\leq
u\leq T:   S_u \leq \frac{1}{\conststime} \mbox{ or } \vola_u \geq \conststime \mbox{
or } \vola_u \leq \frac{1}{\conststime}\}.
\end{eqnarray*}

Consider the following BSDE:
\begin{eqnarray} Y_{t} &= &
H^{N,\conststime}\!-\!\int_{t}^{\stime_\conststime}X_{s}dS_{s}+ \lambda
\int_{t}^{\stime_\conststime}X_{s}^{2}d M_{s}- \sum_{i} \int_{t}^{\stime_\conststime}
\chi_{i,s} dG_{i,s}\label{BSDEeq}
\end{eqnarray} for $0 \leq t \leq \stime_\conststime$ in which
 $H^{N,\conststime} =\mathbf{E}\left( H_T^N | \mathcal{F}_{\stime_\conststime} \right)$. It can be re-written as
\begin{eqnarray} \label{BSDEeqZ}
H^{N,\conststime} & = & Y_{t}-\lambda\int_{t}^{\stime_\conststime} Z_{1,u}^{2} \BSDEdrift_u du +
 \sum_{i=1}^3 \int_{t}^{\stime_\conststime} Z_{i,u} d B_{i,u}\end{eqnarray}
 with
\begin{eqnarray}\label{ZXchi2}Z_{i,u} &=&\volac_{i} \vola_u S_u X_{u} - \xic_i \Mvol(M_t)
X_{u}^2+\sum_{j=1,2} \psi_{i,j,t} \chi_{j,u}
 \end{eqnarray} for $i=2,3$, $Z_{1,u} = \volac_{1}
\vola_u S_u X_{u}$ and $\BSDEdrift_u = \frac{\Mdrift(M_u)}{\volac_1^2
\vola_u^2 S_u^2}$, in which $$\Mdrift(x) = \smallconstant \gamma (\funM^{-1}(x) + \eta) \funM'(\funM^{-1}(x)) + \frac{1}{2} \smallconstant \funM''(\funM^{-1}(x))  \xif(\funM^{-1}(x))^2$$ and $\Mvol(x) =\smallconstant  \xif(\funM^{-1}(x))^2\funM'(\funM^{-1}(x)).$ Note that the change of variable from
$(X,\chi_{1},\chi_2)$ to $(Z_{1},Z_{2},Z_3)$ is one-to-one because
$\psi_t$ is invertible. Since $\frac{ \Mdrift(M_u)}{\vola_u^2
S_u^2}$ is bounded on $[0,\stime_\conststime]$ and $H^{N,\conststime}\in
\Linfty(\mathcal{F}_{\stime_\conststime})$, there exists a pair
$(Z,Y)_{0\leq t\leq{\stime_\conststime}}$ of predictable processes satisfying
BSDE \ref{BSDEeqZ}
 by Theorem 2 of Briand and Hu \cite{BH2006}. Extend these
 processes to $[0,T]$ by setting $Y_t = Y_{\stime_\conststime}$ and $Z_t = 0$
 for $t \geq \stime_\conststime.$

 Define $X$ and $\chi$ in terms
$Z$ with Equation \ref{ZXchi2}. For $m \geq 0$, define
$\overline{X}^m = X 1_{\{|X| \leq m
 \}}$ and similarly for $\overline{\chi}^m.$ Furthermore, let
 $\overline{Z}^m$ be given by Equation \ref{ZXchi2} with $X$
 and $\chi$ replaced by $\overline{X}^m$ and
 $\overline{\chi}^m.$ By Lemma \ref{lemmacontFVapprox}, there
 exists a sequence $\{(\overline{X}^{m,n}$, $\overline{\chi}^{m,n})\}_n$
 of bounded continuous processes with finite variation converging
 to ($\overline{X}^m$, $\overline{\chi}^m)$ in
 $\Htwo$. Define $\overline{Z}^{m,n}$ in terms of $(\overline{X}^{m,n}$,
 $\overline{\chi}^{m,n})$, then $\overline{Z}^{m,n} \rightarrow \overline{Z}^m$ in $\Htwo$ as $n \to \infty$.
Since $\int \overline{Z}^{m,n} d B \to \int \overline{Z}^m d B$,
we also have that $\int_t^{\stime_\conststime} | Z^{m,n}_s |^2 d s \to
\int_t^{\stime_\conststime} | Z^{m}_s |^2 d s$ in $\Lone$.
 Letting
\begin{eqnarray*}
\overline{Y}^{m,n}_{\stime_\conststime} =  Y_{0}-\lambda\int_{0}^{\stime_\conststime}
(\overline{Z}^{m,n}_{1,u})^{2} \BSDEdrift_u du +
 \sum_{i=1}^3 \int_{0}^{\stime_\conststime} \overline{Z}^{m,n}_{i,u} d B_{i,u} \mbox{ and}\\ \overline{Y}^{m}_{\stime_\conststime} = Y_{0}-\lambda\int_{0}^{\stime_\conststime} (\overline{Z}^m_{1,u})^{2} \BSDEdrift_u du +
 \sum_{i=1}^3 \int_{0}^{\stime_\conststime} \overline{Z}^m_{i,u} d
 B_{i,u},\end{eqnarray*}we find $\overline{Y}^{m,n}_{\stime_\conststime} \to \overline{Y}^{m}_{\stime_\conststime} $ in
 $\Lone$ as $n \to \infty.$ Furthermore,
 $(\overline{X}^{m,n},\overline{\chi}^{m,n},\overline{Y}^{m,n})$
 is a s.f.t.s. since it satisfies Equation \ref{BSDEeq} and
 $[\overline{X}^{m,n},\overline{X}^{m,n}]=[\overline{\chi}^{m,n},\overline{\chi}^{m,n}]=0.$

Since $\overline{Y}^{m}_{\stime_\conststime} \to H^{N,\conststime}$ as $m \to \infty,$ we can
find  a sequence
$(X^{n,\conststime,N},\chi^{n,\conststime,N},Y^{n,\conststime,N})_{n \geq 1}$ of s.f.t.s. for each $\conststime$ and $N$ such
that $Y_{\stime_\conststime}^{n,\conststime,N} \rightarrow H^{N,\conststime} = \mathbf{E}\left( H_T^N |
\mathcal{F}_{\stime_\conststime} \right)$ in $\Lone$.

Since $\mathbf{E}\left( H_T^N | \mathcal{F}_{\stime_\conststime} \right)
\rightarrow H_T^N$ as $\conststime \rightarrow \infty$ a.s. by martingale
convergence,
we also have convergence in $\Lone$ by the Dominated
Convergence Theorem. Finally since $H_T^N$ converges to $x
h(\widetilde{S}_{T}^x)$ when $N$ goes to infinity, we can easily
find a s.f.t.s. sequence $(X^{n},\chi^n,Y^{n})_{n \geq 1}$ such
that $Y_T^{n} \rightarrow x h(\widetilde{S}_{T}^x)$ in
$\Lone$ as $n\to \infty.$
\end{proof}

The economic interpretation of Theorem \ref{thcomplete} is that the availability of variance swaps for trading makes the market approximately complete in the sense that any contingent claim with a Lipschitz payoff function can be approximately replicated.

\section{Analytical properties of the approximate solutions}\label{PropertiesSection}

In the
presence of price impacts, the replicating cost of $x$ units of a
contingent claim is not in general $x$ times the replicating cost
of $1$ unit.
When $h$ be a Lipschitz continuous function, recall that for each $x$
an approximating s.f.t.s. for the approximate replication of $x
h(\widetilde{S}_{T}^x)$ is obtained from the solution of BSDE \ref{BSDEeq}, which we denote by $(X^{x}, \chi^{x},
Y^{x})$ to emphasize the dependence on $x$, with the terminal condition
$\mathbf{E} \left(x h^N(\widetilde{S}_{T}^x)\Big|
\mathcal{F}_{\stime_\conststime} \right)$  for $N$ and $\conststime$ large. The theorems in this section give
analytical properties of these approximate solutions for fixed $\conststime$ and $N$. To alleviate the notation, we omit the $\conststime$'s and
$N$'s in all the expressions in this section (e.g. $\tau = \tau_\conststime, h=h^N$, etc ...) when there is no possible confusion.
 For each $t\leq \stime$ and each $x \in \R$, we define $H_{t}(x)
= \frac{1}{x} Y^{x}_{t}$ as the replicating cost per unit for $x$ units of the
claim with payoff function $h$. Furthermore, we let
$H_{t}(0)= \lim_{x\rightarrow 0} H_{t}(x).$ The next theorem states
that this limit exists and is given by the solution of the replication
problem without trade impacts and liquidity costs of Section \ref{replicasubsec}.
Recall that $(\hat{X},\hat{\chi},\hat{Y})_{0 \leq t \leq T}$
denotes the solution of the BSDE \ref{linearBSDE} with terminal
condition $h(S_{T})$.

\begin{theorem}\label{H0thm} $H_{t}(0)= \hat{Y}_{t} =
\mathbf{E}( h(S_T) |\mathcal{F}_{t})$ and
$\frac{1}{x} X^x \rightarrow \hat{X}$ in $\Ltwo(d
\q \times dt)$ as $x \to 0.$
\end{theorem}

\begin{proof} For each $x$, we let $(Z^x,Y^x)_{0\leq t\leq{\stime_\conststime}}$ be the solution of
BSDE \ref{BSDEeqZ} with terminal condition $\mathbf{E} \left(x h(\widetilde{S}_{T}^x)\Big|
\mathcal{F}_{\stime} \right).$ Using the notation of the proof of Theorem \ref{thcomplete}, we have that
$\BSDEdrift_{u}$ is bounded on $[0,\stime]$, which means there exists a
constant $C>0$ such that $ \BSDEdrift_{u}
(Z_{1,u}^x)^{2}  \leq 
C|Z_{u}^{x}|^{2}.$ Take
$|x|<\frac{1}{4 \lambda C C_N }$. First note that since
$\Big|\!\Big|\mathbf{E} \left(
h(\widetilde{S}_{T}^x)\Big| \mathcal{F}_{\stime} \right)
\Big|\!\Big|_\infty \leq C_N$ we know by the maximum principle
(see \cite{K2000}, Proposition 2.1)  that $|Y_{s}^{x}|\leq\
|x| C_N \leq \frac{1}{4 \lambda C}$ for all $0\leq
s\leq\stime.$ Let $H^x =\mathbf{E} \left(
h(\widetilde{S}_{T}^x)\Big| \mathcal{F}_{\stime} \right)$. In
the proof of Theorem \ref{thcomplete} we have shown that
\begin{eqnarray*}
x H^x = Y^{x}_{t}-\lambda
\int_{t}^{\stime}\BSDEdrift_{u}(Z_{1,u}^{x})^{2}du+\int_{t}^{\stime}
Z^{x}_{u}dB_{u},
\end{eqnarray*} thus
\begin{eqnarray*}
x^{2}(H^x)^2 & = & (Y_{t}^{x})^{2}-2
\int_{t}^{\stime}\left(\lambda
\BSDEdrift_{u}(Z_{1,u}^{x})^{2}Y_{u}^{x}
- \frac{1}{2}|Z_{u}^{x}|^{2}\right) du+2 \int_{t}^{\stime}Y_{u}^{x} Z_{u}^{x} dB_{u}\\
 & \geq & (Y_{t}^{x})^{2}+\int_{t}^{\stime}(1-2 \lambda C Y_{u}^{x})|Z_{u}^{x}|^{2}du+2 \int_{t}^{\stime}Y_{u}^{x}Z_{u}^{x}dB_{u}\\
 & \geq & (Y_{t}^{x})^{2}+\int_{t}^{\stime}\frac{1}{2} |Z_{u}^{x}|^{2}du+2 \int_{t}^{\stime}Y_{u}^{x}Z_{u}^{x}dB_{u}.\end{eqnarray*}
 We have  that 
$\mathbf{E}(\int_{t}^{\stime}\BSDEdrift_{u}(Z_{1,u}^{x})^{2}
d u |\mathcal{F}_{t})$ \begin{eqnarray*}&\leq&
\mathbf{E}(\int_{t}^{\stime} C |Z_{u}^{x}|^{2}du|\mathcal{F}_{t})
\leq 2 C \mathbf{E}(x^{2}
(H^x)^{2}|\mathcal{F}_{t}) \leq 2 C x^{2}
C_N^{2}\end{eqnarray*} by taking expectations.
 Since $
Y_{t}^{x} = \mathbf{E} \left( x
H^x + \int_{t}^{\stime}\BSDEdrift_{u}(Z_{1,u}^{x})^{2}du \Big|
\mathcal{F}_t\right),
$ we find \begin{eqnarray}\label{limitHepsilon}
|\frac{1}{x}Y_{t}^{x}-\mathbf{E}(H^x|\mathcal{F}_{t})|
& \leq & x 2 C C_N^{2}.\end{eqnarray} Since
$h(\widetilde{S}_T^x) \rightarrow h(S_T)$ a.s. as $x
\rightarrow 0,$ we have that
$\mathbf{E}(H^x|\mathcal{F}_{t}) =
\mathbf{E}(h(\widetilde{S}_T^x)|\mathcal{F}_{t})$
converges to  $\mathbf{E}(h(S_T)|\mathcal{F}_{t})$
a.s. as $x \to 0$ by the Dominated Convergence Theorem.
Letting $x$ go to zero in Equation \ref{limitHepsilon}, we
have $H_{t}(0)= \hat{Y}_t = \mathbf{E}( h(S_T)
|\mathcal{F}_{t})$.

For the second part of the theorem, let $(\hat{Z},\hat{Y})_{0\leq t\leq T}$
be the solution of
\begin{eqnarray*}
\hat{Y}_t &= & h(S_T)  -
\sum_{j=1}^3 \int_t^T \hat{Z}_{j,s} d B_{j,s} \quad (0\leq t\leq T).
\end{eqnarray*} Then $\hat{Z}_{j,s} = \left(\volac_j \vola_s
\hat{X}_s S_s + \hat{\chi}_{1,s} \psi_{1,j,s}+\hat{\chi}_{2,s} \psi_{2,j,s} \right)$
for $j=1,2$ and $\hat{Z}_{1,s} =\volac_1 \vola_s \hat{X}_s S_s$. Moreover, $ \mathbf{E} \int_{0}^\stime
|\frac{1}{x} Z_u^x - \hat{Z}_u|^2 d u $
\begin{eqnarray*}
&=& \mathbf{E} |\mathbf{E}( h(S_T)
|\mathcal{F}_{\stime})-H^x|^2 -
\left(\mathbf{E} h(S_T)
- \frac{1}{x} Y^x_{0}\right)^2+ 2
\lambda \mathbf{E} \int_{0}^\stime \BSDEdrift_u
\frac{1}{x^2} (Z_{1,u}^x)^2
(Y_u^x- x \hat{Y}_u) d u \\
&\leq & \mathbf{E}
|h(S_T)-h(\widetilde{S}_T^x)|^2 + \frac{ 4 \lambda C C_N
}{x} \mathbf{E} \int_0^\stime |Z^{x}_u|^2
d u
\end{eqnarray*} which goes to $0$ as $x \rightarrow 0.$
Recall that $|X^x_t| =
\left|\frac{Z^x_{1,t}}{\volac_1 \vola_t S_t} \right| \leq
\frac{1}{\volac_1 \conststime^2}\left|Z^x_{1,t}\right|$  on
$[0,\stime]$ and that the same inequality holds for $\hat{X}_t$ and
$\hat{Z}_{1,t}$ for any $0 \leq t \leq \stime.$ Thus  we
find that $\frac{1}{x} X^x$ converges to $\hat{X}$ in
$\Ltwo(d \q \times d t).$
\end{proof}

The next proposition gives an estimate of the error introduced by
using $\widetilde{S}^{x}$ instead of $S^0.$
\begin{proposition}\label{errorapprox}If $h$ is Lipschitz continuous then
$
\mathbf{E} \left|S^{0}_{T+}
-\widetilde{S}^{x}_{T} \right|^2 = O(x^3)$  as $x \rightarrow 0$.
 In particular, $
\mathbf{E} \left| h(S^{0}_{T+}) -
h(\widetilde{S}^{x}_{T}) \right|^2 = O(x^3)$  as $x \rightarrow 0$.
\end{proposition}
\begin{proof}
In terms of $Z^x$, the process $S^0$ can be decomposed as
\begin{eqnarray*}
S_{T+}^{0} = S_T + \int_{0}^\stime \frac{\Mdrift(M_s)}{\volac_1
\vola_s S_s} Z^x_{1,s} d s + \sum_{i} \int_{0}^\stime
\frac{\xic_i \Mvol(M_s)}{\volac_1 \vola_s S_s} Z^x_{1,s} d
B_{i,s},
\end{eqnarray*} since $Z^x_s=0$ for $s$ outside $[{0},\stime]$, whereas
\begin{eqnarray*}
\widetilde{S}_T^{x} = S_T +  \int_{0}^\stime
\frac{\Mdrift(M_s)}{\volac_1 \vola_s S_s} x \hat{Z}_{1,s} d s +
\sum_{i} \int_{0}^\stime \frac{\xic_i \Mvol(M_s)}{\volac_1 \vola_s
S_s} x \hat{Z}_{1,s} d B_{i,s}.
\end{eqnarray*}

In the proof of the previous theorem, we found
\begin{eqnarray*}
 \mathbf{E} \int_{0}^\stime |\frac{1}{x}
Z_u^x - \hat{Z}_u|^2 d u &\leq & \mathbf{E}
|h(S_T)-h(\widetilde{S}_T^x)|^2 + \frac{ 4 \lambda C C_N
}{x} \mathbf{E} \int_{0}^\stime
|Z^{x}_u|^2
d u \\
&\leq & 2 \lambda x^2 \mathbf{E} \left|
\int_{0}^\stime \hat{X}_u d M_u\right|^2 + 8 \lambda C C_N^3  x =
O(x)
\end{eqnarray*} as $x \rightarrow 0$.  Then, for some
positive constant $\hat{C}$,
\begin{eqnarray*}
\mathbf{E} \left|S^{0}_{T+} -
\widetilde{S}^{x}_{T} \right|^2 &\leq & \hat{C}
\mathbf{E}
\int_{0}^\stime \left| x \hat{Z}_u- Z_u^x \right|^2 d u \\
&\leq & x^2 f(x)
\end{eqnarray*} in which $f(x) = O(x)$  as $x \rightarrow
0$.
\end{proof}

Under the additional assumption that $h$ is differentiable we have
that $H_{t}(x)$ is also differentiable at $x=0$ and its derivative can be computed in terms
of the solution of the replication problem without trade impacts. The interpretation of $H'_{t}(0)$ is analogous to the liquidity premium per share $M_{t}$ of the stock.
It gives the additional cost per unit for the replication of the contingent claim due to illiquidity when the number of units replicated is small, i.e. $H_t(x) \approxeq H_t(0) + H'_{t}(0) x$ when $x$ is small. This is comparable to the price of the stock per share $S_t(x) = S_t(0) + M_t x.$

\begin{proposition} Let $0 \leq t \leq \stime$. If $h$ is differentiable everywhere
 except at a finite number of points, then $H_{t}(x)$ is a.s.
differentiable at $x=0$ and
\begin{eqnarray*} H'_{t}(0)&=& \lambda \mathbf{E}\left( \int_{t}^{\stime}
\Mdrift(M_s) \hat{X}_s^2 d s \Big| \mathcal{F}_{t} \right)
-2 \lambda \mathbf{E} \left( h'(S_T)\mathbf{1}_{\{S_T
\leq N\}} (\int_{t}^{\stime} \hat{X}_s d M_s) \Big| \mathcal{F}_{t}
\right).
\end{eqnarray*}
\end{proposition}

\begin{proof}
For $x >0$ small enough, we have that

$ \Big|\frac{1}{x}\left(\frac{Y_{t}^x}{x} -
\hat{Y}_{t} \right) - \lambda \mathbf{E}\left(
\int_{t}^{\stime} \Mdrift(M_s) \hat{X}_s^2 d s \Big| \mathcal{F}_{t}
\right)$

$\quad + 2 \lambda \mathbf{E} \left( h'(S_{T})
\mathbf{1}_{\{S_T \leq N\}} (\int_{t}^{\stime} \hat{X}_s d M_s) \Big|
\mathcal{F}_{t} \right) \Big|$
\begin{eqnarray*}
& = &
\Big|\frac{1}{x}\left(\frac{Y_{t}^x}{x} -
\hat{Y}_{t} \right) - \lambda \mathbf{E}\left(
\int_{t}^{\stime} \BSDEdrift_s \hat{Z}_{1,s}^2 d
s \Big| \mathcal{F}_{t} \right)  + 2 \lambda \mathbf{E} \Big( h'(S_{T})
\mathbf{1}_{\{S_T \leq N\}} (\widetilde{S}_{T}^x-S_T) \Big| \mathcal{F}_{t} \Big) \Big| \\
& \leq & \frac{2 \lambda}{x} \mathbf{E} \left(
\left| h^N(\widetilde{S}_{T}^x) - h^N(S_{T}) -
h'(S_{T})\mathbf{1}_{\{S_T \leq N\}}
(\widetilde{S}_{T}^x-S_T)
\right| \Big| \mathcal{F}_{t} \right) \\
&& + \lambda \left| \mathbf{E} \left(\int_{t}^{\stime}
\BSDEdrift_s \left(\frac{Z_s^x}{x}\right)^2 d s \Big|
\mathcal{F}_{t} \right) - \mathbf{E}
\left(\int_{t}^{\stime} \BSDEdrift_s
\hat{Z}_s^2 d s \Big| \mathcal{F}_{t} \right) \right| \\
&\leq &  \frac{2 \lambda }{x} \mathbf{E} \left(
\left| h^N(\widetilde{S}_{T}^x) - h^N(S_{T}) -
h'(S_{T})\mathbf{1}_{\{S_T \leq N\}}
(\widetilde{S}_{T}^x-S_T) \right| \Big| \mathcal{F}_{t}
\right) \\ &&  + \lambda \left| \mathbf{E}
\left(\int_{t}^{\stime} \BSDEdrift_s
\left(\frac{(Z_s^x)^2}{x^2} - \hat{Z}_s^2\right) d s \Big|
\mathcal{F}_{t} \right) \right|.
\end{eqnarray*} We know that the second term in the last
expression goes to zero when $x \rightarrow 0.$ On the
other hand,
\begin{eqnarray*}
\lim_{x \rightarrow 0} \frac{1}{x}\left(
h^N(\widetilde{S}_{T}^x) - h^N(S_{T}) \right) & = &
\lim_{x \rightarrow 0} \frac{1}{x}\left( h^N(S_{T}+
x \int_{t}^\stime \hat{X}_s d M_s) - h^N(S_{T})\right)\\&=&
h'(S_T) \mathbf{1}_{\{S_T \leq N\}} \int_{t}^\stime \hat{X}_s d M_s
\mbox{ a.s. }
\end{eqnarray*} since $h^N$ is differentiable everywhere except at a
finite number of points. Furthermore, note that
\begin{eqnarray*}
\frac{1}{x}\left| h^N(\widetilde{S}_{T}^x) -
h^N(S_{T}) \right| \leq \hat{C} \left|\int_{t}^\stime \hat{X}_s d M_s
\right|
\end{eqnarray*} in which $\hat{C}$ is the Lipschitz constant of
$h$. We then get the result by the Dominated Convergence Theorem.
\end{proof}

\section{Conclusion}\label{conclu}

This paper extends the liquidity risk model of Çetin et al. \cite{CJP2004} by
hypothesizing the existence of a supply curve that evolves randomly in time and by
studying the impact of trades on the supply curve. This leads to a new
characterization of self-financing trading strategies and a
sufficient condition for no arbitrage. We show the direct connection between stochastic volatility and illiquidity. As a result, contingent claims whose payoffs depend on the
value of the asset can be approximately replicated with the use
of variance swaps. The
replicating costs of such payoffs are obtained from the solutions
of BSDEs with quadratic growth. We show that the marginal cost and the liquidity premium of contingent claims
can be easily computed from the solution of the replication problem without trade impacts.

\section{Acknowledgements}
This work was done at the Center for Applied Mathematics at Cornell University and
is part of the author's Ph.D. thesis under the guidance of Philip
Protter and Robert Jarrow. The author is also very thankful to the two anonymous referees and the associate editor for their constructive comments, and to
 the participants of the Cornell Liquidity Risk Conference and the First European Summer School in Financial Mathematics, and the Bachelier seminar where the results of this paper were presented.

\end{document}